\input amssym.def
\input amssym.tex

\font\teneufm=eufm10
\font\seveneufm=eufm7
\font\fiveeufm=eufm5
\newfam\eufmfam
\textfont\eufmfam=\teneufm
\scriptfont\eufmfam=\seveneufm
\scriptscriptfont\eufmfam=\fiveeufm
\def\goth#1{{\fam\eufmfam\relax#1}}
\font\teneusm=eusm10
\font\seveneusm=eusm7
\font\fiveeusm=eusm5
\newfam\eusmfam
\textfont\eusmfam=\teneusm
\scriptfont\eusmfam=\seveneusm
\scriptscriptfont\eusmfam=\fiveeusm
\def\scr#1{{\fam\eusmfam\relax#1}}
\vskip 3 cm
\magnification = \magstep1

\centerline {\bf SPHERICAL STEIN SPACES}
\bigskip
\bigskip
\centerline {\sl Dmitri Akhiezer and Peter Heinzner}
\footnote{}{Research supported by SFB/TR12 
``Symmetrien und Universalit\"at
in mesoskopischen Systemen'' of the Deutsche
Forschungsgemeinschaft. 
First author supported in part by the Russian Foundation for Basic Research,
Grants 01-01-00709, 02-01-01041. Second author supported in part by
the Schwerpunkt `` Globale Methoden
in der komplexen Geometrie'' of the DFG.}
\bigskip
\bigskip
\centerline {\bf Abstract}
\medskip 
Let $X$ be an irreducible reduced complex space
on which a connected compact Lie group $K$ acts by holomorphic
automorphisms.
Let $G$ be the complexification of $K$
and $\goth g$ the Lie algebra of $G$.
Following the theory of algebraic transformation groups,
we call the complex space $X$ 
spherical if $X$ is normal and its tangent space at some point
is generated 
by the vector fields from a Borel subalgebra $\goth b \subset \goth g$. 
We give several characterizations of spherical Stein spaces.
In particular, we prove that a connected Stein manifold $X$
is spherical if and only if the algebra of $K$-invariant
differential
operators on $X$ is commutative.
\bigskip
\bigskip
\centerline{\bf 1. Introduction}
\medskip
Let $X=(X,\scr O)$
be a complex space, $T_x(X)$ the tangent space at $x\in X$,
and ${\scr T} = {\scr T}_X$ the tangent sheaf of $X$.
We will consider
an action of a connected
compact Lie
group $K$ on $X$. 
We tacitly assume 
that our action  
is continuous
and that each element of $K$
acts as a holomorphic transformation of $X$.
The
complexification of $K$ will be denoted by $G$. By definition,
$G$ is a reductive algebraic group over $\Bbb C$ containing
$K$ as a maximal compact subgroup. 
The action of $K$ gives rise to
a local holomorphic action of $G$ 
on $X$, and so we have the associated Lie homomorphism 
${\goth g} \to {\scr T}(X)$.
Though this homomorphism need not be injective
even if $K$ acts effectively,
we sometimes regard the elements of 
$\goth g$ as vector fields on $X$.
   
Suppose now that $Y$ is an irreducible algebraic variety 
with an algebraic $G$-action and 
$L$ an algebraic line bundle on $Y$ 
with $G$-linearization. Assume that
a Borel subgroup $B\subset G$ has an open
orbit on
$Y$.
Then it is known that $\Gamma (Y,L)$ 
is a multiplicity free $G$-module,
i.e., each irreducible $G$-module occurs in $\Gamma (Y,L)$
with multiplicity 0 or 1. 
This was first shown in [Se], and the argument is as follows.
If $B$ has two weight vectors $s_1, s_2\in \Gamma (Y,L)$
of the same weight then 
$f = s_1/s_2$ is a $B$-invariant rational function
which should be constant on the open $B$-orbit and
hence everywhere on $Y$. Thus a highest weight vector in
$\Gamma (Y,L)$ is uniquely, up to a scalar multiple,
determined by its weight. 

In particular,
the algebra of regular functions ${\Bbb C}[Y]$
is a multiplicity free $G$-module. Moreover,
for $Y$ affine the converse is also true, see [VK].
Our first goal is a generalization of these results to complex spaces.

We recall that an irreducible algebraic $G$-variety $Y$ is called spherical
if $Y$ is normal and $B$ acts on $Y$ with an open orbit.
Similarly, we say that an irreducible reduced complex space
$X$ is spherical (under a $K$-action) if $X$ is normal and
\smallskip
(a) {\it there exists a point $x\in X$ such that
$T_x(X)$
is generated by the elements of $\goth b
$},
\smallskip
\noindent
where $\goth b$ is a Borel subalgebra of $\goth g$.
Since any two Borel subalgebras are conjugated by
an element of $K$, (a) is independent of the choice of $\goth b$. 
\bigskip
{\bf Theorem 1}\ {\it 
Let $X$ be an irreducible reduced complex space,
$K$ a connected compact Lie group acting on $X$, and $L$
a holomorphic line bundle on $X$ with $K$-linearization.
Assume that
$X$ satisfies {\rm (a)}. Then 
${\Gamma } (X,L)$ is a multiplicity free $K$-module. In particular,
\smallskip
{\rm (b)} \ $ \scr O (X)$ is a multiplicity free
$K$-module.}
\bigskip
\noindent
We are mainly interested in $K$-actions on Stein spaces.
For a reduced Stein space $X$  
one has 
the  
complexification theorem, see [He]. Namely,   
there exists another reduced Stein space
$X^{\Bbb C}$, on which
$G = K^{\Bbb C}$
acts
holomorphically, and a $K$-equivariant open
embedding $i: X \hookrightarrow X^{\Bbb C}$, such that 
$i(X)$ is a Runge domain in $X^{\Bbb C}$
and $G\cdot i(X) = X^{\Bbb C}$. 

\bigskip
{\bf Theorem 2}\ {\it Let $X$ be an irreducible reduced Stein space.
Then 
{\rm (a)} and {\rm (b)} are equivalent and 
each of these conditions is equivalent to
\smallskip
{\rm (c)} \ $X^{\Bbb C}$ is an affine algebraic variety
and $B$ has an open orbit on $X^{\Bbb C}$.
\smallskip
\noindent
If 
$X$ is normal then {\rm (a)}, {\rm (b)} and {\rm (c)}
are equivalent to 
\smallskip
{\rm (d)} $X$ is a $K$-invariant domain
in a
spherical affine $G$-variety.}
\bigskip
\noindent
The proof of
(b)$\Rightarrow $(c) for Stein manifolds is sketched in [HW], see p.266 - 267.
We refer the reader to [Ma] for the definition (due to A.Grothendieck)
of a linear differential operator
with holomorphic coefficients on a complex space $X$.
We recall that for local models of complex spaces that definition
is equivalent to the following one. Let $(X, {\scr O}_X)$
be a model space in a Stein domain $U\subset {\Bbb C}^n$,
defined by a coherent ideal sheaf ${\scr J}\subset {\scr O}_U$.
Then linear differential operators with holomorphic coefficients on $X$
are such operators (in the usual sense) on $U$,
which preserve $\scr J$, modulo those which map ${\scr O}_U$
in ${\scr J}$.

Let ${\scr D}_k = {\scr D}_{k,X}$
and ${\scr D} = {\scr D}_X $
be the sheaves of germs of such
differential operators of order not exceeding $k$ and,
respectively, of any order.
For a $K$-sheaf
$\scr F$ on $X$ and
for a $K$-stable subset $U\subset X$
we denote by ${\scr F}(U) ^K$
the set of all $K$-invariant sections of $\scr F$ on $U$.
The sheaves ${\scr D }_k $ and ${\scr D}$
are $K$-sheaves. Moreover,    
${\scr D} (X)$ is a ${\Bbb C}$-algebra on which $K$ acts 
as a group of automorphisms.
In what follows, we consider ${\scr T }(X)$ as a $K$-submodule
of ${\scr D} (X)$.
Our second goal is a characterization of spherical Stein manifolds
in terms of invariant differential operators. 
Again, such a characterization is known for non-singular affine algebraic
varieties, see [HU],[Ag].

\bigskip
{\bf Theorem 3}\ {\rm (i)}{\it If $X$ 
is a (not necessarily reduced) Stein $K$-space satisfying {\rm (b)}
then the algebra of invariant differential operators 
${\scr D}(X)^K$
is commutative. In particular, 
${\scr D}(X)^K$
is commutative for spherical Stein spaces. 

{\rm (ii)}If $X$ is a connected Stein manifold
and ${\scr D}(X)^K$ is commutative then $X$ is spherical.}
\bigskip
\noindent
In the smooth case, (i) is not knew. 
Moreover,
instead of compactness of $K$, one may assume that
the union of irreducible invariant Hilbert
subspaces in ${\scr O}(X)$ is total, see [FT], p.390.
Assertion (ii)
answers a question raised by J.-J. Loeb
in a conversation with the first author several years ago.

Let ${\scr D}_{alg}(Y)$ be the algebra
of algebraic linear differential operators
on an algebraic variety $Y$. Suppose
that $Y$ is acted on by a reductive group $G$ and consider
the subalgebra of $G$-invariant elements ${\scr D}_{alg}(Y)^G
\subset
{\scr D}_{alg}(Y)$. For $Y$ spherical, not necessarily 
affine, ${\scr D}_{alg}(Y)^G$ is a polynomial algebra in
$r$ generators, where $r$ is the rank of $Y$, see [Kn]. 
It is therefore interesting to ask whether 
${\scr D}(X)^K$ is commutative for all spherical complex $K$-spaces.
On the other hand, even if $Y$ is affine, 
the commutativity of ${\scr D}_{alg}(Y)^G$ does not imply 
sphericity of $Y$ for varieties with
(normal) singularities.
Moreover, it can happen
that all
$G$-orbits on $Y$ have positive codimension. The example is actually
contained in [BGG] and is explained below,
see Sect.4, example 3. 
This shows
that (ii) in Theorem 3 is no longer true for Stein spaces
with singularities. 

The authors would like to thank Friedrich Knop for calling their attention
to the paper [Ag] and for useful remarks.
\bigskip
\bigskip
\centerline{\bf 2. Preliminary results}
\medskip
Let $z_1,\ldots,z_n$ be the coordinates in ${\Bbb C}^n$.
We write 
$ \partial _j $ instead of ${\partial \over \partial z_j}$ and
set
$\vert \alpha \vert = \alpha _1 + \ldots +\alpha _n,
\quad \partial
^\alpha = \partial _1 ^{\alpha _1}\cdot \ldots \cdot 
\partial _n ^{\alpha _n}$
for any
$\alpha = (\alpha _1,\ldots, \alpha _n) \in {\Bbb Z}_+^n $.

\bigskip

{\bf Lemma 1}\ {\it 
Let $\Omega \subset {\Bbb C}^n$  be a domain containing the origin
and $E_1,\ldots,E_n$ holomorphic vector fields in $\Omega $. Assume that
$$E_i = \partial _i + F_i,\quad F_i(0) = 0 \quad (i=1,\ldots,n). $$
Then for any $k\ge 1$ and
for any sequence $j =\{j_l\}_{l=1}^k$ of length $k$
with $j_l\in\{1,\ldots,n\}$
one has
$$E_{j_1}\cdot \ldots \cdot E_{j_k} = \partial _{j_1}\cdot \ldots
\cdot \partial _{j_k} \ + \ \sum _{\vert \alpha \vert = k}\ c_{j,\alpha}
 \partial^\alpha  \ + \ G_j, \eqno{(*)}$$
where $c_{j,\alpha}\in {\scr O}(\Omega ),\ c_{j,\alpha}(0) = 0$
and $G_j \in {\scr D}_{k-1}(\Omega )$.}
\bigskip
\noindent
{\it Proof}\ \ Multiplying $(*)$ by $E_i$ on the left we get
$$E_i\cdot E_{j_1}\cdot \ldots \cdot E_{j_k}
= \partial _i\cdot \partial_{j_1}\cdot \ldots \cdot \partial_{j_k}
+ F_i\cdot \partial_{j_1}\cdot \ldots \cdot \partial_{j_k}
+ \sum _{\vert \alpha \vert = k}\ c_{j,\alpha} E_i\cdot \partial ^\alpha +$$
$$
+ \bigl \{ \sum _{\vert \alpha \vert = k}\ (E_i c_{j,\alpha})
\partial ^\alpha + E_i\cdot G_j\bigr \},$$
where the operator in brackets has order $\le k$.
Since $F_i(0) = 0$ and $c_{j,\alpha}(0) = 0$,
our assertion follows by induction.
\hfill $\square $
\bigskip

{\bf Lemma 2}\
{\it With the notation above, let $f\in {\scr O}(\Omega )$ be such a function
that $E_i f  \in {\scr O}(\Omega )f$  for all $i,\ i=1,\ldots, n$. 
If, in addition, $f(0)=0$ then $f=0$.}
\bigskip
\noindent
{\it Proof}\ \
It is easily seen by induction that
$$(E_{j_1}\cdot \ldots \cdot E_{j_k})f \in {\scr O}(\Omega )f$$
for all sequences $j_1, \ldots, j_k$. 
We will show that all derivatives of $f$ vanish at the origin.
For $k \ge 1$ assume by induction
that $\partial ^\alpha f(0) = 0$ for all $\alpha $ with $\vert \alpha \vert
< k$. Then, by Lemma 1,
$$
(\partial _{j_1}\cdot \ldots \cdot \partial_{j_k})f\, (0) =
(E_{j_1}\cdot \ldots \cdot E_{j_k})f\, (0)
 = 0,$$
and so we obtain
$$\partial ^\alpha f(0) = 0$$
for all $\alpha $ with $\vert \alpha \vert = k$. 
\hfill $\square $
\bigskip
{\bf Lemma 3} {\it Let $X$ be a connected
Stein manifold, $K$ a compact connected Lie group
acting on $X$, and
$f\in {\scr O}(X)^K$ a non-constant function.
Then there exists $D\in
{\scr T}(X)^K $, such that $Df\ne 0$.}

\bigskip 
\noindent
{\it Proof}\ \ We identify $X$ with its image $i(X)$ in $ X^{\Bbb C}$.
Note that $X^{\Bbb C}$ is also non-singular.
Furthermore, any $K$-invariant holomorphic function on $X$
extends to a $G$-invariant holomorphic function on $X^{\Bbb C}$,
see [He]. For this reason, we may assume that $X = X^{\Bbb C}$.

Let $X//G$ be the categorical quotient and $\pi: X \to X//G$
the quotient map. 
Take a closed orbit $G\cdot x \subset X$, 
put $p = \pi(x)$ and write $H$ for the
isotropy subgroup $G_x$.
Then $H$ is reductive by Matsushima-Onishchik theorem and
the action of $G$ is locally described by the
holomorphic version of Luna's slice
theorem,
see [Sn] or [He]. Namely,
there is a locally closed $H$-stable Stein subspace $S \subset X$
containing $x$, such that
 the natural map of the twisted product $G\times _H S \to
 X$
is biholomorphic onto a $\pi$-saturated open Stein subset $U\subset X$. 
Since $X$ is non-singular,
there is an $H$-module 
$V$ and an $H$-equivariant biholomorphic map
$\varphi: S \to \varphi (S) \subset V$
of $S$ 
onto an $H$-stable domain $\varphi (S)$ containing the origin.
We will identify $S$ with $\varphi (S)$ and $U$ with $G\times _H S$.

Let $z_1,\ldots,z_n$ be the linear coordinate functions on $V$ and
$E =
\sum_{j=1}^n \ z_j{\partial  \over \partial z_j}$
the Euler vector field. If
$q_d  \in {\Bbb C}[V]^H$ is homogeneous of degree $d>0$
then
$E q_d = d q_d$
by Euler's theorem. But $E $ gives rise to a $
G$-invariant vector field on $U$
tangent to the fibers of the twisted product $G\times _H S$.
Call this vector field again $E$. Then it follows that
$Ef \ne 0$ for any non-constant $G$-invariant holomorphic function on $U$.
This proves our result for $U$ in place of $X$.    

For any coherent 
$G$-sheaf $\scr F$ on $X$ the sheaf $(\pi _*{\scr F})^G$ is coherent on $X//G$,
see [HH].
In particular, 
this applies to $(\pi _*{\scr T})^G$ and 
(by the definition of the complex structure on $X//G$) to
$(\pi _*{\scr O})^G$. 
Now, $\pi (U)$ is an open neighborhood of $p\in X//G$ which can be
taken arbitrarily small. In particular,
we may assume that $\pi (U)$ has 
Runge property for any coherent
sheaf on the Stein space $X//G$.
The assertion of the lemma follows 
then from Runge approximation theorem
for coherent sheaves, see [GR], Kap. V, $\S$ 6. Indeed, the subspaces 
${\scr O}(X)^G = (\pi _*{\scr O})^G (\pi (X)) 
\subset (\pi _*{\scr O})^G (\pi (U)) = {\scr O}(U)^G$
and
${\scr T}(X)^G = (\pi _*{\scr T})^G(\pi (X)) 
\subset (\pi _*{\scr T})^G (\pi (U)) = {\scr T}(U)^G$
are dense 
and the Lie derivative $(D,f)\to Df$ is continuous 
in both arguments  
in the canonical Fr\'echet topology.
\hfill $\square $
\eject
\bigskip
\bigskip
\centerline{\bf 3. Proofs of the theorems}
\medskip 
\noindent
{\it Proof of Theorem 1}\ 
Take a point $x \in X$,
such that $T_x(X)$ is generated by global vector fields
from $\goth b$. In a neighborhood of $x$ we have a 
non-singular analytic set
through $x$ whose tangent space
coincides with $T_x(X)$. It follows that $x$
is a non-singular point. Choose a coordinate
neighborhood $\Omega $ of $x$, denote by
$z_1, \ldots, z_n$ the local coordinates satisfying
$z_i(x) = 0$, and write $\partial _i $ for $\partial \over \partial z_i$
in $\Omega $.  
Then there exist $E_1, \ldots , E_n\in  \goth b$
such that $E_i (x)=\partial_i$ for all $i$.

Let ${\scr L}$ denote the sheaf of germs
of holomorphic sections of $L$.
Without loss of generality
assume that ${\scr L}\vert \Omega  \simeq {\scr O}$.
Let $s_0 \in {\scr L}(\Omega )$ be a non-vanishing section.
 
Note that ${\scr L}(X)$ and ${\scr O} (X)$ are
$(K,\goth g)$-modules.
Furthermore, the action of $\goth g$ is local in the sense
that for any open set 
$U\subset X$ the algebra $\goth g$ acts on ${\scr O}(U)$ and ${\scr L}(U)$
commuting with restriction maps. For the multiplication map
${\scr O}(U)\times {\scr L}(U) \to {\scr L}(U)$ and any element of $\goth g$
we have the Leibniz rule.

Assume that ${\scr L}(X)$ 
contains two isomorphic irreducible $K$-submodules which
do not coincide as subspaces. Then these modules have
the same highest weight with respect to $\goth b$. 
Let $s_1, s_2\in {\scr L} (X)$ be the corresponding weight vectors.
The sections 
$s_1,s_2 $
are linearly independent, and one has
$$Es_k = \lambda (E)s_k \quad (k=1,2) $$
for all $E\in {\goth b}$, where $\lambda :{\goth b}\to {\Bbb C}$
is a Lie algebra character. 
Take a non-trivial linear combination
$s = c_1s_1 + c_2s_2$,
such that $s(x) =0$.
Restricting $s$ to $\Omega $, we can write
$s\vert {\Omega } = fs_0$, where $f\in {\scr O}(\Omega )$.
Since $Es_0 = \varphi _E s_0$, where $\varphi _E \in {\scr O}(\Omega )$,
we get 
$$ Ef = (\lambda (E) - \varphi _E)f\in {\scr O}(\Omega )f$$
by Leibniz rule.
But then Lemma 2 shows that $f = 0$,
contradictory to the fact that
$s_1,s_2$ are linearly independent.
\hfill $\square $

\bigskip
\noindent
{\it Remark}\ \ Here is another proof of Theorem 1
for Stein spaces.  
Condition (a) implies in this case that $B$ has an open orbit on $X^{\Bbb C}$.
In particular, ${\scr O}(X^{\Bbb C})^G = {\Bbb C}$.
Since any two closed $G$-orbits are separated by invariant functions,
$G$ has exactly one closed orbit on $X^{\Bbb C}$. It
follows that $X^{\Bbb C}$ is ($G$-equivariantly biholomorphic
to) an affine algebraic variety on which $G$ acts algebraically,
see [Sn], Cor. 5.6.
On the other hand, the line bundle $L$ extends to
a holomorphic line bundle $\tilde L$
on $X^{\Bbb C}$ with $G$-linearization
and, moreover, the multiplicities of irreducible $K$-modules in
$\Gamma (X,L)$ and $\Gamma (X^{\Bbb C}, \tilde L)$ are the same,
see [HH], Cor. 3 and Identity Principle.   
The assertion follows now from the corresponding
result for affine algebraic varieties.
\hfill $\square $
\bigskip
\noindent
{\it Proof of Theorem 2}\ \ ${\rm (a)}\Rightarrow {\rm (b)} $ is already
proved.
\smallskip
${\rm (b)}\Rightarrow {\rm (c)} $.\ Since 
${\scr O}(X^{\Bbb C})$
is a (dense) $K$-invariant subspace in ${\scr O}(X)$,
it is a multiplicity free $K$-module and,
equivalently,  
a multiplicity free $G$-module.  
In particular, ${\scr O}(X^{\Bbb C})^G ={\Bbb C}$.
As in the above remark, we see that $X^{\Bbb C}$
is an affine algebraic variety with algebraic action of $G$.
The algebra of regular functions
on this variety is a multiplicity free $G$-module. Therefore $B$ acts on
$X^{\Bbb C}$
with an open orbit, see [VK].
\smallskip
${\rm (c)}\Rightarrow {\rm (a)} $
and
${\rm (d)}\Rightarrow {\rm (a)} $ 
are obvious.
\smallskip
${\rm (c)}\Rightarrow {\rm (d)} $ follows from the fact each point
of $X^{\Bbb C}$ is contained 
in the $G$-orbit of some point of $i(X)$. Therefore,
if $X$ is normal then $X^{\Bbb C}$ is also normal.
\hfill $\square $

\bigskip
\noindent
{\it Proof of Theorem 3}\ (i) Let $\hat K$ denote the set
of equivalence classes of irreducible linear representations of $K$.
For each $\delta \in \hat K$
we have the isotypic component ${\scr O}_\delta (X) 
\subset {\scr O}(X)$ of type $\delta $.
By a theorem of Harish-Chandra, ${\scr O}_\delta (X)$ is closed in the 
canonical Fr\'echet topology of ${\scr O}(X)$ and the sum
of all subspaces ${\scr O}_\delta (X)$ is dense, see [Ha-Ch]
or [Akh], \S5.1.
  
Clearly, each 
${\scr O}
_\delta (X)$ is stable under ${\scr D}(X)^K$.
Condition (b) means that 
${\scr O}_\delta (X)$ is irreducible.
Therefore any $D\in {\scr D}(X)^K$
acts on ${\scr O}_\delta (X)$ as a scalar operator
by Schur's lemma.
It follows that if
$D$ is of the form $D =[D_1,D_2]$, where $D_1,D_2 \in {\scr D}(X)^K$,
then $Df = 0$ for all $f\in {\scr O}(X)$.
For a Stein space, this implies $D=0$. Indeed, 
let $x\in X$ be any point,
$D_x$ the germ of $D$ at $x$ and 
$g_x\in {\scr O}_x$ any element of the local ring.
In some neighborhood $U$ of $x$ one can represent 
$g_x$ by a function
$g\in {\scr O}(U)$.
Choose $U$ to be a Runge domain in $X$ so that
the image of the restriction map ${\scr O}(X)
\to {\scr O}(U)$ is dense in the canonical Fr\'echet
topology.
Since $D$ is continuous in this topology   
we get $Dg = 0$ showing that $D_x = 0$.

(ii) Suppose $X$ is a connected Stein manifold with 
commutative algebra ${\scr D}(X)^K$.
In order to prove that $X$ is spherical
it suffices to show
that $X^{\Bbb C}$ is an affine algebraic $G$-variety. Indeed,
since
$X^{\Bbb C}$ is obviously non-singular, the 
commutativity of ${\scr D}_{alg}(X^{\Bbb C})$
implies then that $X^{\Bbb C}$ is spherical, see [Ag], Satz 2.5.
On the other hand,
we know that ${\scr O}(X^{\Bbb C})^G ={\Bbb C}$
implies that 
$X^{\Bbb C}$ is an affine algebraic $G$-variety
(see the above remark). 
Therefore,
since ${\scr O}(X^{\Bbb C})^G = {\scr O}(X)^K $,  
we only have to check that ${\scr O}(X)^K = {\Bbb C}$.  

Assume the contrary. Let $f\in {\scr O}(X)^K$ 
be a non-constant function. By Lemma 3 there is
an invariant holomorphic vector field $D$
on $X$, such that $Df \ne 0$. Denote by $M_f$ the operator (of order 0)
of multiplication by $f$. Then $M_f$ is also invariant and
$[D,M_f] = M_{Df} \ne 0$. This contradicts
the commutativity of ${\scr D}(X)^K$.
\hfill $\square $   
\bigskip
\bigskip
\centerline{\bf 4. Some examples}
\bigskip
1)\ Let $X$ be a compact complex manifold without non-constant
meromorphic functions, ${\rm dim}\, X \ge 2$,
and let $K=\{\rm id\}$ be the trivial transformation group of $X$.
For any holomorphic line bundle $L$ 
the space $\Gamma (X,L)$ has dimension 0 or 1. Thus
$\Gamma (X,L)$ is a multiplicity free $K$-module.
However, $X$ is obviously non-spherical,
showing that the converse statement to Theorem 1 is false. 
Due to Rosenlicht's theorem, this
phenomenon does not occur in the algebraic setting.
More precisely, suppose that a non-singular algebraic variety
$X$ is non-spherical with respect to
an algebraic action of a reductive group $G$. Then
there exist a line bundle $L$ on $X$ and a $G$-linearization of $L$,
such that some irreducible $G$-module occurs in $\Gamma (X,L)$
with multiplicity greater than 1.

\bigskip
2)\ Let $\Gamma $ be a discrete Zariski dense subgroup in 
$G =K^{\Bbb C}$ and $X = G/\Gamma $. Then it is known that
${\scr O}(X) = \Bbb C$, see e.g.
[Akh], \S 5.3. In this situation, one has 
a natural isomorphism between ${\scr D}(X)$
and the universal enveloping algebra $U(\goth g)$.
Indeed, any $D\in {\scr D}(X)$ can be viewed as
a holomorphic differential operator on $G$ commuting
with right translations $R_g$ for $g\in \Gamma$.
For any $f\in {\scr O}(G)$ define a function $\varphi = \varphi _f
\in {\scr O}(G\times G)$ by
$$\varphi (g,x) = (DR_{g^{-1}}f)(x) - (Df)(xg^{-1}).$$
Then for $\gamma \in \Gamma $ we have 
$$\varphi (g\gamma, x) = 
(DR_{\gamma ^{-1}}R_{g^{-1}}f)(x) - (Df)(x\gamma ^{-1}g^{-1}) 
= $$
$$=
R_{\gamma ^{-1}}(DR_{g^{-1}}f - R_{g^{-1}}Df)(x) =
\varphi (g,x\gamma ^{-1}).$$
Since ${\scr O}(X) = {\Bbb C}$,
it follows that
$$\varphi (g,x) = \varphi (e,xg^{-1}) = 0,$$
showing that $D$ commutes with all right translations.

The isomorphism ${\scr D}(X) = U(\goth g)$ is compatible
with the $G$-action. Thus
${\scr D}(X)^K $ is isomorphic to the center of $U(\goth g)$.
Though ${\scr D}(X)^K$ is commutative, $X$ 
is spherical 
only if $K$ is a torus.
Thus, in Theorem 3 (ii) it is essential that
$X$ is Stein.
\bigskip  
3)\ Following [BGG], consider the normal surface
$$Y =\{z=(z_1,z_2,z_3) \in {\Bbb C}^3
\ \vert \  z_1^3+z_2^3+z_3^3 = 0\}$$
with $G ={\Bbb C}^*$ 
and $K = \{\zeta \in {\Bbb C}^* \, \vert \,
\vert \zeta \vert = 1\}$
acting on $Y$ by $z \to \zeta\cdot z$.
The associated infinitesimal operator is the Euler vector field
$E = z_1\partial _1 + z_2\partial _2 + z_3\partial _3 \in {\scr D}(Y)$.
For any $p\in {\Bbb Z},\ p\ge 0,$ set
$${\scr O}^{(p)} (Y) = \{f\in {\scr O}(Y) \, \vert
\, Ef = pf\},\quad 
{\scr D}^{(p)} (Y) = \{D\in {\scr D}(Y) \, \vert \, [E,D] = pD\}.$$
Note that ${\scr O}^{(p)}(Y)$ is the space
of all homogeneous functions of degree $p$ on $Y$,
in particular, ${\scr O}^{(p)}(Y) \subset {\Bbb C}[Y]$.
On the other hand, it is easily seen that
$${\scr D}^{(p)}(Y) = \{ D\in {\scr D}(Y) \, \vert \, D\cdot 
{\scr O}^{(q)}(Y)
\subset {\scr O}^{(p+q)}(Y) \quad {\rm for \ all } \ q\ge 0 \},$$
hence
${\scr D}^{(p)}(Y)\subset {\scr D}_{alg}(Y)$. 
A differential operator $D \in {\scr D}(Y)$ is $G$-invariant
if and only if $[E,D] = 0$. 
Thus
$${\scr D}(Y)^K = {\scr D}_{alg}(Y)^G = 
 {\scr D}^{(0)}(Y).$$
One of the results in [BGG]
yields an isomorphism of ${\Bbb C}$-algebras
${\scr D}_{alg}(Y)^G \simeq {\Bbb C}[E]$.
This shows that assertion (ii) in Theorem 3
is false for singular Stein spaces. 

For the same reason, its algebraic 
counterpart is false if an affine variety 
has (normal) singularities.   
Moreover,
one can say why
the proof in the non-singular case
does not go through for singular varieties
(cf. [Bi], Lemma 4.2). 
That proof is based on the density property of ${\scr D}_{alg}(Y)$.
Namely, given $l$ linearly independent
functions $f_1, \ldots, f_l \in {\Bbb C}[Y]$
and $l$ arbitrary functions $g_1, \ldots, g_l \in {\Bbb C}[Y]$,
there exists an algebraic differential operator $D$
such that $Df_i = g_i$ for all $i,\ i= 1,\ldots, l.$
In general, this is no longer true
for varieties with singularities, e.g., for the surface $Y$
considered above.  
Indeed, it is shown in [BGG]
that any $D \in {\scr D}_{alg}(Y)$ can be written as 
a finite sum $D = D_0+D_1+D_2+\ldots $,
where 
$D_p\in {\scr D}^{(p)}(y)$,
i.e., $D_p$ raises the homogeneity degree of a function by $p$.
Thus density fails already for $l=1$.
\eject
\bigskip
\bigskip
\centerline {\bf References}
\medskip
\noindent
[Ag] I.Agricola, {\it Invariante Differentialoperatoren und die
Frobenius-Zerlegung einer $G$-Variet\"at}, J. of Lie Theory {\bf 11}
(2001), 81 - 109.
\medskip 
\noindent
[Akh] D.N.Akhiezer, {\it Lie group actions in complex analysis},
Vieweg, Braunschweig - Wiesbaden, 1995.
\medskip
\noindent
[BGG] I.N.Bernstein, I.M.Gel'fand, S.I.Gel'fand, {\it Differential
operators on a cubic cone}, Uspekhi Mat. Nauk {\bf 27} (1972),
no.1, 185 -
190 (Russian); English translation:
Russian Math. Surveys, {\bf 27} (1972), no.1,
169 -174.
\medskip
\noindent
[Bi] F.Bien, {\it Orbits, multiplicities and differential operators},
Contemp. Math., Amer. Math. Soc. vol. {\bf 145}, 1993, 199 - 227.
\medskip
\noindent
[FT] J.Faraut, E.G.F.Thomas, {\it Invariant Hilbert spaces of holomorphic
functions}, J. of Lie Theory {\bf 9} (1999), 383 - 402
\medskip
\noindent
[GR]\ H.Grauert, R.Remmert,
{\it Theorie der Steinschen R\"aume},
Springer-Verlag, Berlin - Heidelberg - New York, 1977.
\medskip
\noindent  
[Ha-Ch] Harish-Chandra, {\it Discrete series for semisimple Lie
groups }II, Acta Math. {\bf 116} (1966), 1 - 111.
\medskip
\noindent
[HH] J.Hausen, P.Heinzner, {\it Actions of compact groups
on coherent sheaves}, Transformation Groups {\bf 4} (1999),
25 - 34.
\medskip 
\noindent
[He] P.Heinzner, {\it Geometric invariant theory on Stein spaces},
Math. Ann. {\bf 289} (1991), 631 - 662.
\medskip
\noindent
[HU] R.Howe, T.Umeda, {\it The Capelli identity, the double commutant
theorem and multipli- \-city-free actions},
Math. Ann. {\bf 290} (1991), 565 - 619.
\medskip
\noindent
[HW] A.T.Huckleberry, T.Wurzbacher, {\it Multiplicity-free complex manifolds},
Math. Ann. {\bf 286} (1990), 261 - 280.
\medskip
\noindent
[Kn] F.Knop, {\it A Harish-Chandra homomorphism for reductive group actions},
Ann. of Math. (2) {\bf 140} (1994), 253 - 288.
\medskip
\noindent
[Ma] B.Malgrange, {\it Analytic spaces}, 
l'Enseignement Math. (2) {\bf 14} (1968), 1 - 28.
\medskip
\noindent
[Se] F.J.Servedio, {\it Prehomogeneous vector spaces and varieties},
Trans. Amer. Math. Soc. {\bf 176} (1973), 421 - 444.
\medskip
\noindent
[Sn] D.M.Snow, {\it Reductive group actions on Stein spaces},
Math. Ann. {\bf 259} (1982), 79 - 97.
\medskip
\noindent
[VK] E.B.Vinberg, B.N.Kimel'feld, {\it Homogeneous domains in flag manifolds
and spherical subgroups of semisimple Lie groups},
Funktsional'nyi Analiz i ego Prilozheniya {\bf 12} (1978), 12 - 19
(Russian); English translation:
Funct. Anal. and Appl. {\bf 12} (1978), 168 - 174.

\noindent
------------------------------------------------------------------------
\bigskip
\noindent
D.Akhiezer: Institute for Information Transmission Problems, B.Karetny 19,

\noindent
101447 Moscow, Russia
\bigskip
\noindent
P.Heinzner: Ruhr-Universit\"at Bochum, Fakult\"at f\"ur Mathematik,
Universit\"atsstra\ss e 150, 

\noindent 44780 Bochum, Germany

\end